\documentclass{elsarticle}

\usepackage{
hyperref}

\usepackage{amssymb, amsmath}
\usepackage[sc]{mathpazo}
\usepackage{bm} 
\usepackage{latexsym} 
\usepackage{graphics} 
\usepackage{url} 
\usepackage{enumerate}

\usepackage[all]{xy}

\usepackage{graphicx}

\usepackage{float}

\usepackage{color}

\usepackage[utf8]{inputenc}

\journal{ }

\newtheorem{theo}{Theorem}

\newtheorem{lemm}[theo]{Lemma}

\newtheorem{exam}[theo]{Example}

\newtheorem{property}[theo]{Property}

\newenvironment{proo}[1][\proofname]{\normalfont{\itshape
#1{:}}\quad\mdseries\ignorespaces}
{{$\Box$}{\vskip\belowdisplayskip}}
\newcommand{\proofname}{Proof}

\newtheorem{defi}[theo]{Definition}

\newtheorem{rem}[theo]{Remark}

\newtheorem{problem}[theo]{Problem}

\DeclareMathOperator{\rank}{rank}

\begin{document}

\begin{frontmatter}

\title{A Rellich's result revisited and sensitivity of solutions of parametrized linear systems}

\author[mymainaddress]{Jos\'e Carlos Bellido}
\ead{josecarlos.bellido@uclm.es}

\author[mysecondaryaddress]{Luis Felipe Prieto-Mart\'inez}
\ead{luisfelipe.prieto@upm.es}

\address[mymainaddress]{Departamento de Matem\'aticas, ETSI Industriales, INEI, Universidad de Castilla-La Mancha, Campus Universitario S/N, E13071, Ciudad Real, Spain}
\address[mysecondaryaddress]{Departamento de Matem\'atica Aplicada, ETS Arquitectura, Universidad Polit\'ecnica de Madrid, Avd. Juan de Herrera 4, E28040, Madrid, Spain}

\begin{abstract}
In this paper we revisit a result due to Franz Rellich on smoothness of solutions of parametrized linear systems. With this result as a starting point, we obtain finer smoothness results in an elementary fashion and propose an efficient adjoint algorithm for computing sensitivities of $(n-1)$-deficient systems,  being $n$ the order of the system.  
\end{abstract}


\end{frontmatter}


\section{Introduction}

More than 50 years ago Franz Rellich pioneered the investigation on the problem of perturbation of eigenvalue problems with respect to matrix system parameters, both for finite and infinite dimensional systems \cite{R}.  In this paper we focus our attention on the following simple lemma from \cite{R} in the finite dimensional framework, where smoothness with respect to a single parameter of solutions of homogeneous finite dimensional linear systems verifying a normalization constraint is established. 

\begin{theo}[Rellich's Lemma] \label{theo.rellich} Let $\bold D(\varepsilon)=[\gamma_{ij}(\varepsilon)]_{1\leq i,j\leq n}$ be a $n\times n$ matrix whose coefficients, $\gamma_{ij}(\varepsilon)$, for $i,j=1,\ldots,n$, are real analytic functions in a neighborhood of some $\varepsilon_0$ and such that for each $\varepsilon$ in this neighborhood, $det(\bold D(\varepsilon))=0$. Then, for some neighborhood of $\varepsilon_0$,  there exist analytic functions $\alpha_1(\varepsilon),\ldots,\alpha_n(\varepsilon)$ such that the column vector $\bold x(\varepsilon)=[\alpha_1(\varepsilon), \ldots  \alpha_n(\varepsilon)]^T$ satisfies:
\begin{enumerate}
\item[a)] \label{theo.rellich_1}$\bold D(\varepsilon)\bold x(\varepsilon)=\bold 0$, and 
\item[b)] \label{theo.rellich_2}$\|\bold x(\varepsilon)\|=1$.
\end{enumerate}
\end{theo}

In this paper we revisit this result, whose proof is elementary, exploring the generalization of this theorem to the multi-parameter case, i.e. when coefficient matrix $\bm D$ depends on a vector of variables $\bm \varepsilon\in \mathbb R^N$, for both for homogeneous and non-homogeneous linear systems. More concretely, we would like to provide answers to the following questions. 

\begin{problem} \label{problem.fefferman} Let $\bm D(\bm\varepsilon)=[\gamma_{ij}(\bm\varepsilon)]_{1\leq i,j\leq n}$,  for $\bm\varepsilon\in U\subset \mathbb R^N$.  If for $1\leq i,j\leq n$, $\gamma_{ij}(\bm\varepsilon)$ is analytic (respectively in $C^l(U)$), decide:

\begin{itemize}
\item Given $\bm b(\bm\varepsilon)=[b_1(\bm\varepsilon),\ldots, b_n(\bm\varepsilon)]^T$, with $b_i(\bm \varepsilon)$ analytic (resp. in $C^l(U)$) for $1\leq i\leq n$, does the system 
\begin{equation} \label{eq.1nh} \bm D(\bm\varepsilon)\bm x(\bm\varepsilon)=\bm b(\bm\varepsilon),\end{equation}
admits analytic solutions?

\item does the linear system
\begin{equation} \label{eq.1}  \bold D(\bm \varepsilon)\bold x(\bm \varepsilon)=\bold 0\end{equation}
 admits a solution $\bm x(\bm\varepsilon)=[x_1(\bm\varepsilon),\ldots, x_n(\bm\varepsilon)]^T$ satisfying
 \begin{equation} \label{eq.2}\|\bold x(\bm \varepsilon)\|=1.\end{equation}
and such that  each $x_i(\bm\varepsilon)$ is analytic (resp. is in $C^1(U)$)?

\end{itemize}

\end{problem}

A second step in the direction set in Problem \ref{problem.fefferman} would be to device an efficient way of computing derivatives or differentials of solutions of linear systems with respect to parameters. Thus, the second question we address in this investigation is the following.

\begin{problem} \label{problem.derivative}  Let $\bm D(\bm\varepsilon)=[\gamma_{ij}(\bm\varepsilon)]_{1\leq i,j\leq n}$,  for $\bm\varepsilon\in U\subset \mathbb R^N$.  If, for $1\leq i,j\leq n$, the coefficients $\gamma_{ij}(\bm\varepsilon)$ are analytic (resp. are in $C^1(U)$) for $1\leq m\leq N$, determine the value of $\frac{\partial \bm x}{\partial \varepsilon_m}(\bm\varepsilon_0)$,  where 
$\bm x(\bm\varepsilon)$ is one of the analytic solutions (resp. $C^1(U)$) of \eqref{eq.1} satisfying \eqref{eq.2}, and such that $x(\bm \varepsilon_0)=\bm u$ (with $\bm u$ a fixed given unitary vector such that $\bm{D}(\bm{\varepsilon})\bm{u}= \bm{0}$).
\end{problem}

Both, studying smoothness of solutions and the effective computation of derivatives, usually referred as sensitivities, fit into the field of \emph{sensitivity analysis}, which for this kind of linear problems goes back to the times of A. Turing, who in his influential paper \cite{T}, pointed out the interest on the problem of sensitivity of solutions of linear systems (in that paper, he also introduced the definition of condition number). The literature on the subject itself, or on intimately related topics, is really vast. Rellich's result that motivates this paper is contextualized as an auxiliary result in the infancy of the perturbation theory of eigenproblems \cite{R}. \cite{K} is a classical reference on the subject. Sensitivity of eigenvalues and eigenvectors of linear problems is a matter of great practical interest in many engineering contexts, as for instance and just for citing one among many, structural design \cite{B1,B2}. Of course, the problems we study could be seen as particular cases of eigenvector sensitivity analysis, although our questions are even more elementary and our aim here is to check how far we can get pushing forward the elementary ideas in the proof of Theorem \ref{theo.rellich}. Literature on this subject is really overwhelming both from the mathematical and engineering sides, and to include here an exhausting list of references is out of the scope of this paper. We additionally reference \cite{Lancaster,Seyranian}. 

Another very related subject is singular value decomposition of parametrized matrices, where both the smoothness of the decomposition and methods for computing it are addressed. The literature on this topic is also huge. We cite \cite{Volker91,Volker93,Sun1,Sun2}. 

Coming back to our very elementary problem on smoothness and sensitivity of solutions of parametrized linear systems, it is worth mentioning the highly cited survey \cite{D}. See also \cite{ROHN}. The standard approach in these works is the following: starting from a system $\bm D\bm x=\bm b$, with $\bm D$ a non-singular matrix, consider the perturbed system $(\bm D+\bm{\Delta D})\bm x=\bm b+\bm{\Delta b}$, where $\bm{\Delta D}$ and $\bm{\Delta b}$ are interpreted as perturbations or errors. Thus, these problems are usually studied from the Numerical Linear Algebra viewpoint and from what  is called \emph{interval analysis} in the literature. This viewpoint is mainly oriented to computing bounds for the uncertainty of the components of $\bm x$ rather than to compute derivatives.

On the contrary, we assume coefficients of $\bm D,\bm b$  are {\it smooth} functions depending on a vector of parameters $\bm\varepsilon\in\mathbb R^N$, instead of considering perturbations. This is the natural framework in many applied situations, where the coefficients of the matrix sometimes depend on some parameters with some (physical, or economical, or engineering) meaning. Furthermore, in the case of homogeneous systems the constraint \eqref{eq.2} is taken into consideration. Smoothness and explicit computation of so-called {\it frame solutions}, i.e. families of orthonormal solutions of the homogeneous system, has been studied intensively before \cite{ColeSoren84,Rheinboldt88,Rheinboldt93}. In this paper, we will provide a simple proof for the existence of smooth frame solutions and an efficient method for computing sensitivities in the multi-parameter case for $(n-1)$-deficient systems (i.e. $\rank(\bm D({\bm\varepsilon}))=n-1$).

The existence of smooth solutions to linear systems with respect to parameters has received some recent attention due to its relationship with \emph{Whitney's Extension Problem}. We highlight the recent articles  \cite{FK, FL1, FL2}, for the case $U=\mathbb R^N$, where a complete characterization of the linear systems admitting such a solution is obtained. We are not so interested in obtaining such characterizations, but simple and coarse criterions for the existence of smooth solutions, following the philosophy of \cite{R}, and to device methods for sensitivity computation. 

Outline of the paper is the following.  In Section \ref{section.proofrellich} we include an updated version of the proof of Theorem \ref{theo.rellich} and, after this,  we discuss the questions in Problem \ref{problem.fefferman} as generalizations of Theorem \ref{theo.rellich}: the necessity of the analyticity condition  (Subsection \ref{section.analyticity}), the multi-parameter case (Subsection \ref{section.multiparameter}) and finally  the general non-homogeneous case (Subsection \ref{section.general}). In Section \ref{section.ff} we study the existence of  smooth frame solutions when $\rank\,({\bm D}({\bm\varepsilon}))<n-1$. Finally, Section \ref{section.final} is devoted to the exposition of two sensitivities computation algorithms for the multi-parameter case and for $(n-1)$-deficient systems. We provide a direct method, inspired in Nelson's Method for simple eigenvector derivative calculation \cite{Nelson}, and an adjoint method (more efficient for large values of $N$) inspired by \cite{B1,B2}.

\section{Rellich's Theorem} \label{section.proofrellich}

In this section, first, we include, for readers' convenience and to understand the new results in this paper, the proof of Theorem \ref{theo.rellich}. In Subsections 2.1, 2..2 and 2.3 extensions and generalizations of Theorem \ref{theo.rellich} are given, partially answering Problem \ref{problem.fefferman}.

\begin{proo}[Proof of Theorem \ref{theo.rellich}] Let us denote by $U$ the neighborhood of $\varepsilon_0$ where the entries of matrix $\bm D(\varepsilon)$ are analytic and such that $\det (\mathbf{D}(\varepsilon))=0$ for all $\varepsilon\in {U}$. Set $r=\max_{\varepsilon\in U}(\rank(\bm D(\varepsilon)))$. We may assume that $\bm D(\varepsilon)$ is not the trivial matrix in $U$.  So $1\leq r\leq n-1$.  

First we construct an analytic solution verifying Theorem \ref{theo.rellich}, 
a). There is no loss in generality (performing a permutation of the equations and of the variables) in assuming that   $\det([\gamma_{ij}(\varepsilon)]_{i,j=1}^r)$ is a minor which is not constantly equal to 0 in $ U$. 

For $1\leq i,j\leq n$,   denote by $\Gamma_{ij}(\varepsilon)$ to  the cofactor  of $\gamma_{ij}(\varepsilon)$ in the submatrix $[\gamma_{ij}(\varepsilon)]_{i,j=1}^{r+1}$. Defining
\begin{equation} \label{eq.f_k}f_k(\varepsilon)=\begin{cases} \Gamma_{r+1,k}(\varepsilon)&\text{for }k=1,\ldots, r+1\\ 0 & \text{for }k=r+2,\ldots,n, \end{cases} 
\end{equation}
functions $f_{k}(\varepsilon)$ are analytic  and not simultaneously  constantly zero (because $f_{r+1}(\varepsilon)=\Gamma_{r+1,r+1}(\varepsilon)\neq 0 $ for some $\varepsilon \in U$). Moreover, we have that for $i=1,\ldots, n$,
$$[\gamma_{i1}(\varepsilon),\ldots,\gamma_{in}(\varepsilon)]   \begin{bmatrix}f_1(\varepsilon)\\ \vdots \\ f_n(\varepsilon) \end{bmatrix}=\sum_{k=1}^{r+1} \gamma_{ik}(\varepsilon)\Gamma_{r+1,k}(\varepsilon)=0,$$
since the second term in the equality is the determinant of the  $(r+1)\times(r+1)$ matrix obtained by replacing the $(r+1)$-st row of $[\gamma_{ij}(\varepsilon)]_{i,j=1}^{r+1}$ by $[\gamma_{i1}(\varepsilon),\ldots, \gamma_{i,r+1}(\varepsilon)]$ and therefore vanishes:
\begin{itemize}
	\item for $1\le i< r$, since it is the determinant of a matrix with two equal rows;
	\item for $i=r$,  since it is the determinant of matrix $[\gamma_{ij}(\varepsilon)]_{i,j=1}^{r+1}$;
	\item for $r<i\le n$, since it is the determinant of a matrix such that its last row is linear combination of the rest of rows.
\end{itemize}

In order to complete the proof let us assume, for the shake of simplicity, that $\varepsilon_0=0$. There exist some $m\in\mathbb{N}\cup\{0\}$ such that, for $1\leq i\leq n$, the power series $f_i(\varepsilon)$ is of order at least $m$. So we have:
\begin{equation} \label{eq.factorcomun} (|f_1(\varepsilon)|^2+\ldots+|f_n(\varepsilon)|^2)=\varepsilon^{2m}(h_0+h_1\varepsilon+\ldots), \qquad h_0\neq 0
 \end{equation}
where the power series $h_0+h_1\varepsilon+\ldots$ converges for sufficiently small $|\varepsilon|$. Consequently:
$$(|f_1(\varepsilon)|^2+\ldots+|f_n(\varepsilon)|^2)^{\frac{1}{2}}=\varepsilon^m(\omega_0+\omega_1\varepsilon+\ldots), \qquad \omega_0\neq 0 $$
where the power series $\omega_0+\omega_1x+\ldots=\sqrt{h_0+h_1x+\ldots}$ converges for small $|\varepsilon|$. Since, as we have established, the power series expansion of each $f_i(\varepsilon)$ are of order at least $m$, we can define 
$$\alpha_i(\varepsilon)=\frac{f_i(\varepsilon)}{\varepsilon^m(\omega_0+\omega_1\varepsilon+\ldots)},\qquad i=1,\ldots, n $$
which is also a convergent power series for small $|\varepsilon|$ and, for some $1\leq i\leq m$, $\alpha_i(\varepsilon)\neq 0$.

\hspace{12cm}\end{proo}

In the second part of the proof, the following elementary property of the set of formal power series has been essential. This property will be recalled later.

\begin{property} \label{property.obstruction} If $f_1(\varepsilon),\ldots, f_k(\varepsilon)$ are convergent power series in a neighborhood of $\varepsilon_0\in\mathbb R$, then there exists  some $m\in \mathbb{N}\cup \{0\}$, such that $ \alpha_1(\varepsilon)=\frac{f_1(\varepsilon)}{(\varepsilon-\varepsilon_0)^m},\ldots,  \alpha_k(\varepsilon)=\frac{f_k(\varepsilon)}{(\varepsilon-\varepsilon_0)^m}$ are all of them convergent  power series and for some $i$, $1\leq i\leq k$, $\alpha_i(\varepsilon_0)\neq 0$.

\end{property}

\medskip
In the following subsections, we address generalizations of Theorem \ref{theo.rellich} in different directions.

\subsection{Analyticity assumption in Rellich's Theorem} \label{section.analyticity}

As was already pointed out by Rellich, and due to the important role played by Property \ref{property.obstruction} in the proof of Theorem \ref{theo.rellich}, it is not possible to replace the analyticity condition in Theorem \ref{theo.rellich} by a weaker smoothness condition. That is, if the coefficients  $\gamma_{ij}(\varepsilon)\in C^l(U)$, $\l\ge 1$, for $1\leq i,j\leq n$, in general it is not possible to find $\alpha_1(\varepsilon),\ldots, \alpha_n(\varepsilon)\in C^l(U^*)$, with $U^*$ being a neighborhood of $\varepsilon_0$ and satisfying the conditions $a)$ and $b)$ in Theorem \ref{theo.rellich}. The following example is an adaptation of the one appearing in \cite{R} (recall that in \cite{R} the eigenproblem is addressed) illustrating this fact.

\begin{exam} \label{ex.rellich} For $n=2$, consider the following matrix, which entries are continuous and have continuous derivatives of all orders in $\mathbb R$:
$$\bm D(\varepsilon)=\begin{cases}  \begin{bmatrix}e^{-\frac{1}{\varepsilon^2}}(1-cos\frac{2}{\varepsilon}) & -e^{-\frac{1}{\varepsilon^2}}sin\frac{2}{\varepsilon}\\
-e^{-\frac{1}{\varepsilon^2}}sin\frac{2}{\varepsilon} & e^{-\frac{1}{\varepsilon^2}}(1+cos\frac{2}{\varepsilon}) \end{bmatrix} & \text{for }\varepsilon\neq 0 \\ \\ \begin{bmatrix} 0 & 0 \\ 0 & 0 \end{bmatrix} & \text{for }\varepsilon=0 \end{cases}$$

\noindent For $\varepsilon\neq 0$ and for some real function $\lambda(\varepsilon)$, any solution takes the form:
$$\bm x(\varepsilon)=\lambda(\varepsilon)\cdot \begin{bmatrix} cos(1/\varepsilon)\\ sin(1/\varepsilon))\end{bmatrix} $$

\noindent  There is no neighborhood $U^*$ of $\varepsilon_0=0$ with a vector $\bm  x(\varepsilon)=\begin{bmatrix}\alpha_1(\varepsilon)\\ \alpha_2(\varepsilon)\end{bmatrix}$ such that $\alpha_1(\varepsilon), \alpha_2(\varepsilon)$ are continuous and such that for every $\varepsilon\in U^*$ they satisfy  a) and b) in Theorem \ref{theo.rellich}.

\end{exam}

What is possible is, starting from a matrix with $C^l(U)$ coefficients, to obtain a vector $\bm x(\varepsilon)$ which components are also in  $C^l(U)$ and such that satisfy $a)$ in Theorem \ref{theo.rellich} and different o zero in some neighborhood of $\varepsilon_0$, but not Condition $b)$. Example \ref{ex.rellich} also illustrates this situation. Starting from a matrix whose coefficients are $C^\infty(\mathbb R)$, there exist a column vector $\bm x(\varepsilon)$ (which is not, in general, unique) satisfying $a)$ and such that its components are also $C^\infty(\mathbb R)$. Take, for instance,
$$\bm x(\varepsilon)=e^{-\frac{1}{\varepsilon^2}}\cdot \begin{bmatrix} cos(1/\varepsilon)\\ sin(1/\varepsilon))\end{bmatrix}. $$

\noindent In other words, it is possible to obtain the following result in the spirit of Theorem \ref{theo.rellich} for weaker smoothness conditions.

\begin{theo} \label{lemmatilde}  Let $\varepsilon_0\in \mathbb R$, and let $U$ be a neighborhood of $\varepsilon_0$.  Let $\bm D(\varepsilon)=[\gamma_{ij}(\varepsilon)]_{1\leq i,j\leq n}$ a matrix such that  each $\gamma_{ij}(\varepsilon)\in C^l(U)$, for  all $i,j=1,\ldots,n$,  and such that $det(\bm D(\varepsilon))=0$ for all $\varepsilon\in U$.  Then:

\begin{enumerate}

\item There exist  functions  $f_1(\varepsilon),\ldots,f_n(\varepsilon)$ in $C^l(U)$ such that the column vector 
$$\bm x(\varepsilon)= \begin{bmatrix}f_1(\varepsilon)\\ \vdots \\ f_n(\varepsilon) \end{bmatrix}$$
satisfies
$$\bm{D}(\varepsilon) \bm x(\varepsilon)=\bm 0$$
for any $\varepsilon \in U$;
\item Let $r=\max_{\varepsilon\in U}(\rank(\bm D(\varepsilon)))$ and assume that  $\rank(\bm D(\varepsilon_0))=r$. Then, there exists a neighborhood $U^*$ of $\varepsilon_0$ such that for $\varepsilon\in U^*$, $\bm x(\varepsilon)\neq \bm 0$, and consequently, there exist functions $\alpha_1(\varepsilon),\ldots,\alpha_n(\varepsilon)$ in $C^l(U^*)$ such that the column vector 
$$\bm x(\varepsilon)= \begin{bmatrix}\alpha_1(\varepsilon)\\ \vdots \\ \alpha_n(\varepsilon) \end{bmatrix}$$
satisfies conditions a) and b) in Theorem \ref{theo.rellich}, for any $\varepsilon\in U^*$.

\end{enumerate}

\end{theo}

\begin{proo} The proof follows the lines of that of Theorem \ref{theo.rellich}. Assume, again, that $0<r<n$.  There is no loss in generality assuming, also,  that   $\det([\gamma_{ij}(\varepsilon)]_{i,j=1}^r)$ is a minor which is non trivial in some neighborhood $U^*\subset U$ of $\varepsilon_0$. The proof of part {(1)} is the same as the one of Theorem  \ref{theo.rellich} (with the obvious modifications) and will not be repeated here. For the proof of Part 2, recovering the previous notation, we just notice that since $\rank(\bm D(\varepsilon_0))=r$, therefore constant and maximal in the whole neighborhood $U^*$, 
$$|f_1(\varepsilon)|^2+\ldots+|f_n(\varepsilon)|^2 $$
is a   $C^l(U^*)$ function, which does not vanish for $\varepsilon\in U^*$. So does
$$(|f_1(\varepsilon)|^2+\ldots+|f_n(\varepsilon)|^2)^{1/2} $$
and so, for $1\leq i \leq n$
$$\alpha_i(\varepsilon)=\frac{f_i(\varepsilon)}{(|f_1(\varepsilon)|^2+\ldots+|f_n(\varepsilon)|^2)^{1/2}}$$
is in $C^l(U^*) $.

\hspace{12cm}\end{proo}

The solutions constructed in the proof fail to have smoothness at the points $\varepsilon_0$ such that $\rank(\bm D(\varepsilon_0))$ is not maximal (is less than $r$). This is the reason of the problems appearing in the smoothness of the eigenvectors when two or more eigenvalues coalesce in the eigenproblem context.

\subsection{Multi-parameter case} \label{section.multiparameter}

Now let us consider the several variables case of the second part of Problem \ref{problem.fefferman}, that is, now $\bm D(\bm\varepsilon)=[\gamma_{ij}(\bm\varepsilon)]_{1\leq i,j\leq n}$ denotes a $n\times n$ matrix such that each entry $\gamma_{ij}$ depends on a vector of variables $\bm\varepsilon\in\mathbb R^N$, $N>1$, and solutions are vectors $\bm x(\bm\varepsilon)$.

The natural analogue of Theorem \ref{theo.rellich} for this case is false, as we can see in the following example. Again, this is a consequence of the fact that analytic functions in several variables do not satisfy an analogue of Property \ref{property.obstruction}.

\begin{exam} \label{ex.mine} For $n=2,N=2$, consider the following matrix, which entries depend on a vector of variables $\bm\varepsilon=(\varepsilon_1,\varepsilon_2)$ and are analytic functions for every $(\varepsilon_1,\varepsilon_2)\in\mathbb R^2$:
$$\bm D(\varepsilon_1,\varepsilon_2)=\begin{bmatrix}2\varepsilon_1\varepsilon_2& \varepsilon_2^2-\varepsilon_1^2\\ 0 & 0  \end{bmatrix}$$

\noindent Any solution of system 
$$\bm D(\varepsilon_1,\varepsilon_1)\bm x(\varepsilon_1,\varepsilon_2)=\bm 0$$ verifies 
$$\bm x(\varepsilon_1,\varepsilon_2)=\lambda(\varepsilon_1,\varepsilon_2)\cdot\begin{bmatrix}\varepsilon_1^2-\varepsilon_2^2\\ 2\varepsilon_1\varepsilon_2 \end{bmatrix}, $$
for some real function $\lambda(\varepsilon_1,\varepsilon_2).$

\noindent Now if we impose $\bm x(\varepsilon_1,\varepsilon_2)$ to satisfy Equation (2), we obtain
$$\|[\varepsilon_1^2-\varepsilon_2^2, 2\varepsilon_1\varepsilon_2]^T\|=\sqrt{(\varepsilon_1^2-\varepsilon_2^2)^2+(2\varepsilon_1\varepsilon_2)^2}=\varepsilon_1^2+\varepsilon_2^2$$

\noindent so the  solution $\bm x(\varepsilon_1,\varepsilon_2)$ such that $\|\bm x(\varepsilon_1,\varepsilon_2)\|=1$ for $(\varepsilon_1,\varepsilon_2)\neq (0,0)$ is (up to sign)
$$\bm x(\varepsilon_1,\varepsilon_2)=\begin{bmatrix}\frac{\varepsilon_1^2-\varepsilon_2^2}{\varepsilon_1^2+\varepsilon_2^2}\\ \frac{2\varepsilon_1\varepsilon_2}{\varepsilon_1^2+\varepsilon_2^2} \end{bmatrix} $$

\noindent It is not possible to extend such a solution to be continuous at $(0,0)$.

\end{exam}

In a similar fashion to what happened in the previous section, for a matrix $\bm D(\bm\varepsilon)$ which entries are smooth (analytic or in $C^l(U)$), it is possible to find a vector which entries are also smooth (analytic or in $C^l(U)$) and such that it satisfies \eqref{eq.1}, but not \eqref{eq.2}, in general, that is, it is possible to prove an analogue to Theorem \ref{lemmatilde} for several variables, but not of Theorem 1.

\begin{theo} \label{lemmatildeseveral}

 Let $\bm D(\bm\varepsilon)=[\gamma_{ij}(\bm\varepsilon)]_{1\leq i,j\leq n}$ such that, for $i,j=1,\ldots,n$, $\gamma_{ij}(\bm\varepsilon)$ is analytic (resp. $C^l(U)$)   in a neighborhood $U$ of some $\bm\varepsilon_0\in\mathbb R^N$ and such that for each $\bm\varepsilon\in U$, $det(\bm D(\bm\varepsilon))=0$. Then:

\begin{enumerate}

\item There exist functions  $f_1(\bm\varepsilon),\ldots,f_n(\bm\varepsilon)$ which are analytic (resp. $C^l(U)$) and such that the column vector 
$$\bm x(\bm\varepsilon)= \begin{bmatrix}f_1(\bm\varepsilon)\\ \vdots \\ f_n(\bm\varepsilon) \end{bmatrix}$$
satisfies \eqref{eq.1};

\item Let $r=\max_{\bm\varepsilon\in U}(\rank(\bm D(\bm\varepsilon)))$.  Assume that $\rank(\bm D(\bm\varepsilon_0))=r$. Then, there exists a neighborhood $U^*$ of $\bm\varepsilon_0$, such that there exists functions $\alpha_1(\bm \varepsilon),\ldots, \alpha_n(\bm\varepsilon)$ which are analytic (resp. $C^l(U^*)$) and the column vector
$$\bm x(\bm\varepsilon)=\begin{bmatrix} \alpha_1(\varepsilon)\\ \vdots\\ \alpha_n(\bm\varepsilon)\end{bmatrix}$$

satisfies conditions \eqref{eq.1} and \eqref{eq.2} for $\bm\varepsilon\in U^*$.

\end{enumerate}

\end{theo}

The proof follows completely the lines to the one of Theorem \ref{lemmatilde} and will not be repeated.

\subsection{Extension to non-homogeneous linear systems} \label{section.general}

Now we deal with the first question raised in Problem \ref{problem.fefferman}. In this case, the corresponding result has a direct and elementary proof and unifies the one-parameter and the multi-parameter cases:

\begin{theo} \label{theo.nohomo} Let $\bm D(\bm\varepsilon)=[\gamma_{ij}(\bm\varepsilon)]_{1\leq i,j\leq n}$, $\bm b=[b_1(\bm\varepsilon),\ldots, b_n(\bm\varepsilon)]^T$ such that their entries are analytic (resp. $C^l(U)$) in a neighborhood $U$ of some $\bm\varepsilon_0\in\mathbb R^N$, $N\geq 1$. If $\rank(\bm D(\bm\varepsilon_0))=\rank(\bm D(\bm\varepsilon_0)\mid \bm b(\bm\varepsilon_0))=r$ and for every $\bm\varepsilon\in U$ $\rank(\bm D(\bm\varepsilon)\mid \bm b(\bm\varepsilon))\leq r$ ($\star$), there exists some neighborhood  $U^*$ of $\bm\varepsilon_0$ such that  there exist functions $g_1(\bm\varepsilon),\ldots, g_m(\bm\varepsilon)$ which are analytic (resp. $C^l(U^*)$) and such that $\bm x(\bm \varepsilon)=[g_1(\bm\varepsilon),\ldots,g_n(\bm\varepsilon)]^T$ satisfies \eqref{eq.1nh} for $\bm\varepsilon\in U^*$.


\end{theo}

\begin{proo} There is no loss in generality assuming that $det([\gamma_{ij}(\bm\varepsilon_0)]_{1\leq i,j\leq r})$ is a non-trivial minor.   There is a neighborhood $U^*$ of $\bm\varepsilon_0$ contained in $U$ such that for every $\bm\varepsilon\in U^*$, $\det([\gamma_{ij}(\bm\varepsilon_0)]_{1\leq i,j\leq r})\neq 0$. Consider the matrix :
$$\bm{\widetilde D}(\bm\varepsilon)=\left[\begin{array}{l l l | l l l}\gamma_{11}(\bm\varepsilon) & \hdots & \gamma_{1r}(\bm\varepsilon) & \gamma_{1,r+1}(\bm\varepsilon) &\hdots & \gamma_{1n}(\bm\varepsilon)\\ 
\vdots & & \vdots & \vdots & & \vdots \\ 
\gamma_{r1}(\bm\varepsilon) & \hdots & \gamma_{rr}(\bm\varepsilon) & \gamma_{r,r+1}(\bm\varepsilon) &\hdots & \gamma_{rn}(\bm\varepsilon)\\ \hline 0& \hdots & 0 &  \\ \vdots & & \vdots & & \bm I_{n-r} \\ 0 & \hdots & 0 &\end{array}\right]$$

\noindent where $\bm I_{n-r}$ denotes the identity matrix of size $(n-r)\times(n-r)$.

The unique solution of $\bm{\widetilde D}(\bm\varepsilon)\bm x(\bm\varepsilon)=\bm b(\bm \varepsilon)$ is a solution of $\bm D(\bm\varepsilon)\bm x(\bm\varepsilon)=\bm b(\bm\varepsilon)$. So the result is a direct consequence of Cramer's rule applied to the system $\bm{\widetilde D}(\bm\varepsilon)\bm x(\bm\varepsilon)=\bm b(\bm\varepsilon)$.

\hspace{12cm}\end{proo}

Condition ($\star$) holds automatically if $\rank(\bm D(\bm\varepsilon))$ is constant in $U$ and equal to $n$.

\section{System of smooth linearly independent solutions to Problem \ref{problem.fefferman}}\label{section.ff}

Let us begin by introducing some notation. Let $U\subset \mathbb R^N$. For us a \emph{vector field} is a map $\bm x:U\to \mathbb R^n$ (indeed, this has already appeared above). We say that it is  analytic (resp. is $C^l(U)$) if so are its components. The following definition clarifies the term {\it frame field}, which appears in the literature with more than one meaning.

\begin{defi}\label{defi.ff} Let $U\subset\mathbb R^N$, a \textbf{frame field}  is a map $F:U\to(\mathbb R^n)^k$
$$\bm\varepsilon=(\varepsilon_1,\ldots,\varepsilon_N)\longmapsto (\bm x_1(\bm\varepsilon),\ldots,\bm x_k(\bm\varepsilon)) $$

\noindent   such that, for every $\bm\varepsilon\in U$,  $\bm x_i(\bm\varepsilon)\cdot \bm x_j(\bm\varepsilon)=\delta_{ij}$ (so $k\leq n$). Note that each $\bm x_i$ is a $n$-dimensional vector field in $U$. We say that it is analytic (resp. in $C^l(U)$) if so are each vector field $\bm x_i$, for $1\leq i\leq k$.

\end{defi}

\subsection{Frame fields of solutions of  homogeneous linear systems} \label{subsection.ffh}

In the context of the second case of Problem \ref{problem.fefferman}, in this subsection we prove, in some cases, the existence of not only a vector field corresponding to a solution, but a frame field $(\bm x_1(\bm\varepsilon),\ldots, \bm x_k(\bm\varepsilon))$ consisting on $k$ solutions.

\medskip

 The proccess used in this section generalizes an analogue construction done, again, in \cite{R} for eigenvalue problems. Let us start with the following simple observation.
\begin{lemm} \label{lemm.B} Let $\bm D$ be a constant matrix of size $n\times n$ and rank $r\geq 1$. Let $k=n-r$ and $\bm x_1,\ldots,\bm x_k $ be any orthonormal basis of solutions of the linear system $\bm D\bm x=\bm 0$. Then the matrix 
$$\bm B=\bm D+\bm x_k \bm x_k^T$$ 

\noindent  of size $n\times n$  has rank $r+1$ and satisfies 
$$\bm B \bm x_i\begin{cases}=\bm 0 & \text{for }i=1,\ldots, k-1\\ \neq \bm 0 & \text{for }i=k \end{cases}$$

\end{lemm}

\begin{proo}  For $i=1,\ldots, k-1$, taking into account that ${\bm x}_1,\dots,{\bm x}_k$ are in the kernel of $\bm D$ and constitute an orthogonal system we get
$$\bm B\bm x_i=(\bm D+\bm x_k \bm x_k^T) \bm x_i=\bm 0  $$

\noindent and analogously, having in mind, this time, that $\bm x_k\bm x_k=1$
$$\bm B \bm x_k=(\bm D+\bm x_k \bm x_k^T)\bm x_k=\bm x_k\neq \bm 0$$

Since $\bm B\bm x=\bm 0$ has at least $k-1$ solutions, $\rank(\bm B)\leq r+1$. To see the equality, we have to check that there is no more than $k-1$ linearly independent solutions of this system. Suppose that we have a vector $\bm x_0$ such that $\{\bm x_0,\ldots,\bm x_{k}\}$ is orthonormal. Then 
$$\bm B\bm x_0=(\bm D+\bm x_k\bm x_k^T) \bm x_0=\bm D\bm x_0$$

\noindent So, if $\bm B\bm x_0=\bm 0$ then $\bm D\bm x_0=\bm 0$ which leads to a contradiction with  $\rank(\bm D)=r$.

\hspace{12cm}\end{proo}

This idea allow us to prove the following theorem, which is an improved statement of Theorem \ref{theo.rellich} (in the one-parameter case). This result does not appear explicitely in \cite{R}.

\begin{theo}   \label{theo.mainanalytic} Let $U$ be a neighborhood of $\varepsilon_0\in\mathbb R$.  Let $\bm D(\varepsilon)=[\gamma_{ij}(\varepsilon)]_{1\leq i,j\leq n}$ be such that  each $\gamma_{ij}(\varepsilon)$, for $1\leq i,j\leq n$, is  analytic in $U$ and such that for each $\varepsilon\in U$, $det(\bm D(\varepsilon))=0$.  Suppose that $r=\max_{\varepsilon\in U}\{\rank(\bm D(\varepsilon))\}$ and let $k=n-r$. Then there exist an analytic frame field  $(\bm x_1(\varepsilon),\ldots, \bm x_k(\varepsilon))$ such that for each $\varepsilon$ in a neighborhood of $\varepsilon_0$, for $1\leq i\leq k$, $\bm D(\varepsilon)\bm x_i(\varepsilon)=\bm 0$.

\end{theo}

\begin{proo}  We proceed by induction.  The case $k=1$ is a trivial consequence of Theorem \ref{theo.rellich}. 

For $k> 1$, using, again, this result we can obtain a vector that will be denoted by $\bm x_k(\varepsilon)$ such that $\bm D(\varepsilon) \bm x_k(\varepsilon)=\bm 0$ and $\|\bm x_k(\varepsilon)\|=1$. For each $\varepsilon\in U$ consider any orthonormal set of solutions $\{\bm y_1(\varepsilon),\ldots, \bm y_{k-1}(\varepsilon),\bm x_k(\varepsilon)\}$ of the system $\bm D(\varepsilon)\bm x=\bm 0$ containing the one above. Now define $\bm B(\varepsilon)=\bm D(\varepsilon)+\bm x_k(\varepsilon)\bm x_k^T(\varepsilon)$. $\bm B(\varepsilon)$ also satisfies the hypothesis of the theorem, with $\rank(\bm B)\leq r+1$. So we can apply the induction hypothesis to obtain a frame field of $k-1$ vectors $(\bm x_1(\varepsilon),\ldots, \bm x_{k-1}(\varepsilon))$ satisfying the requirements in the theorem. And since, for each $\varepsilon\in U$ they belong to the subspace spanned by $\bm y_1(\varepsilon),\ldots, \bm y_{k-1}(\varepsilon)$, the set
$\{\bm x_1(\varepsilon),\ldots, \bm x_{k-1}(\varepsilon),\bm x_k(\varepsilon)\}$ is an orthonormal set.

\hspace{12cm}\end{proo}

In a similar fashion, we can obtain a similar result replacing analyticity by $C^l (U)$ smoothness in the one-parameter case ($N=1$) and in the multi-parameter case ($N>1$).

\begin{theo}   \label{theo.main}  Let $\bm\varepsilon_0\in \mathbb R^N$, and let $U$ be a neighborhood of this $\bm\varepsilon_0$.  Let $\bm D(\bm\varepsilon)=[\gamma_{ij}(\bm\varepsilon)]_{1\leq i,j\leq n}$ be such that  each $\gamma_{ij}(\bm\varepsilon)$, for $1\leq i,j\leq n$, is in $C^l(U)$ and such that for each $\bm\varepsilon\in U$, $det(\bm D(\bm\varepsilon))=0$.  Suppose that $r=\max_{\varepsilon\in U}\{\rank(\bm D(\bm\varepsilon))\}$ and let $k=n-r$. Then there exist $k$  vector fields  $\bm x_1(\bm\varepsilon),\ldots, \bm x_k(\bm\varepsilon)$ which are $C^l(U)$ such that for each $\bm\varepsilon\in U$, for each $1\leq i,j\leq k$, $i\neq j$,  $\bm x_i(\bm\varepsilon)^T\bm x_j(\bm\varepsilon)=0$ and $\bm D(\bm\varepsilon)\bm x_i(\bm\varepsilon)=\bm 0$.  Moreover, if $\rank(\bm D(\bm\varepsilon_0))=r$, then there exists a neighborhood $U^*$ of $\bm\varepsilon_0$  such that, for $\bm\varepsilon\in U^*$, it is possible to choose these vectors in such a way that $(\bm x_1(\bm\varepsilon),\ldots, \bm x_k(\bm\varepsilon))$ is a frame field.

\end{theo}

The proof follows the lines of the previous one.

\medskip

Let us remark that, although  Theorems \ref{theo.mainanalytic} and \ref{theo.main} are stated as purely existence result, the method exhibited in the proof of Theorem \ref{theo.mainanalytic} is constructive, as shown in the following example, that we hope it helps to clarify the algorithm. Anyway, to compute the frame field using this method is not  computationally efficient, since it requires to compute too many determinants.

\begin{exam} \label{exam.parte1} In the case $n=3, N=1$, consider the matrix $\bm D(\varepsilon)=\begin{bmatrix}  \varepsilon & 0 & 0 \\ 0 & 0 & 0 \\ 0 & 0 & 0 \end{bmatrix}$. 

For $U=\mathbb R$, this matrix satisfies 
$$\rank(\bm D(\varepsilon))=\begin{cases}1 & \text{if }\varepsilon\neq 0 \\ 0 & \text{if }\varepsilon=0\end{cases}$$

\noindent so we expect to obtain a frame field of solutions consisting in 2 vector, defined in $\mathbb R\setminus \{0\}$.

\noindent The proof of Theorem \ref{theo.rellich} explains how to obtain one of the solutions, using certain cofactors. This solution is 
$$\bm x_1(\bm\varepsilon)=\frac{1}{\|[0,\varepsilon,0]^T\|}\begin{bmatrix}0 \\ \varepsilon\\ 0 \end{bmatrix}=\begin{bmatrix} 0 \\ 1 \\ 0 \end{bmatrix}.$$

\noindent Now, following the idea in Theorem 13 we repeat the same proccess, but for the matrix 
$$B(\bm\varepsilon)=\bm D(\varepsilon)+\bm x_1(\varepsilon)\bm x_1(\varepsilon)^T=\begin{bmatrix} \varepsilon & 0 & 0 \\ 0 & 1 & 0 \\ 0 & 0 & 0 \end{bmatrix}$$

\noindent obtaining the solution 
$$\bm x_2(\bm\varepsilon)=\begin{scriptsize}\frac{1}{\|\left[\left|\begin{array}{l l} 0 & 0 \\ 1 & 0\end{array}\right|, \left|\begin{array}{l l} \varepsilon & 0 \\ 0 & 0\end{array}\right|,\left|\begin{array}{l l} \varepsilon & 0 \\ 0 & 1\end{array}\right| \right]^T\|}\begin{bmatrix}\left|\begin{array}{l l} 0 & 0 \\ 1 & 0\end{array}\right|\\ \\ \left|\begin{array}{l l} \varepsilon & 0 \\ 0 & 0\end{array}\right|\\ \\\left|\begin{array}{l l} \varepsilon & 0 \\ 0 & 1\end{array}\right| \end{bmatrix}\end{scriptsize}=\begin{bmatrix} 0 \\ 0 \\ 1 \end{bmatrix}.$$

Note that, in this case, the frame field can be extended to the whole open set $U$, although this may not happen in general.

\end{exam}

Finally, note that  this analytic frame field $(\bm x_1(\bm \varepsilon),\ldots,\bm x_k(\bm \varepsilon))$ obtained with this method, is not the only possible one satisfying the conditions.

\begin{rem}\label{rem.movimiento}  Using Gram-Schmidt Method it is possible to find analytic (resp. in $C^l(U)$) vector fields $\bm{\widetilde x_{k+1}}(\bm\varepsilon),\ldots, \bm{\widetilde x_{n}}(\bm\varepsilon)$ in such a way that 
\begin{equation} \label{eq.orthonormalbasis} \{\bm x_1(\bm \varepsilon),\ldots,\bm x_k(\bm \varepsilon),\bm{\widetilde  x_{k+1}}(\bm\varepsilon),\ldots, \bm{\widetilde x_{n}}(\bm\varepsilon)\}\end{equation}

\noindent is an orthonormal basis.

Let $(\bm y_1(\bm\varepsilon),\ldots, \bm y_k(\bm \varepsilon))$ be a frame field. Then it satisfies the conditions if and only if, there exists a matrix $\bm K(\bm \varepsilon)$ that  maps the elements in the basis $(\bm x_1(\bm \varepsilon),\ldots,\bm x_k(\bm \varepsilon))$ into elements in the basis  $(\bm y_1(\bm\varepsilon),\ldots, \bm y_k(\bm \varepsilon))$ of the form
\begin{equation} \label{eq.K}\bm K(\bm\varepsilon)=\bm P(\bm\varepsilon)\left[\begin{array}{l | l}\bm T(\bm \varepsilon) & \bm 0\\ \hline \bm 0 & \bm I_{n-k}  \end{array}\right]\bm P(\bm \varepsilon)^{-1} \end{equation}

\noindent where $\bm P(\bm \varepsilon)$ is the matrix whose columns are the vectors in the base \eqref{eq.orthonormalbasis}, $\bm T(\bm\varepsilon)$ is a matrix  which entries are analytic (resp. are in $C^l(U)$) and such that for each $\bm\varepsilon\in U$ belongs to the orthogonal group $\mathcal O(k)$ and finally $\bm I_{n-k}$ is the $(n-k)\times (n-k)$ identity matrix. Note that the entries of $\bm y_i(\bm\varepsilon)$ are analytic (resp. are in $C^l(U)$).

A matrix whose entries depend on some variables and  belongs to $\mathcal O(n)$ for each value of these variables, such as $\bm K(\bm\varepsilon)$,  is sometimes called \emph{kinematic matrix}. In fact, we can view the frame field $(\bm y_1(\bm\varepsilon),\ldots, \bm y_k(\bm \varepsilon))$ as the result of applying a rigid motion to the original frame field $(\bm x_1(\bm \varepsilon),\ldots,\bm x_k(\bm \varepsilon))$ in such a way that the entries in the matrix of this rigid motion have the adequate smoothness with respect to the variables in $\bm\varepsilon$.

\end{rem}

\subsection{Linearly independent solutions for non-homogeneous linear systems} 

From  Theorem \ref{theo.nohomo} and the results in the previous subsection, it is easy to prove the following:

\begin{theo} \label{theo.nohomo1}   Let $\varepsilon_0\in \mathbb R$, and let $U$ be a neighborhood of this $\varepsilon_0$.  Let $\bm D(\varepsilon)=[\gamma_{ij}(\varepsilon)]_{1\leq i,j\leq n}$, $\bm b(\varepsilon)=(b_1(\varepsilon),\ldots, b_n(\varepsilon))^T$ be such that, for $1\leq i,j\leq n$, $\gamma_{ij}(\varepsilon)$, $b_i(\varepsilon)$, are  analytic for $1\leq i,j\leq n$. Suppose that $r=\rank(\bm D(\varepsilon_0))=\rank([\bm D(\varepsilon_0)\mid \bm b(\varepsilon_0)])$ and that for every $\varepsilon\in U$,  $\rank([\bm D(\varepsilon)\mid \bm b(\varepsilon)])\leq r$. Let $k=n-r$. Then there exist a neighborhood $U^*$ of $\varepsilon_0$, an analytic vector field $\bm x_p(\varepsilon)$ and an analytic frame field $(\bm x_1(\varepsilon),\ldots, \bm x_k(\varepsilon))$ such that any analytic vector field $\bm x(\varepsilon)$ satisfying Equation \eqref{eq.1nh} can be writen as:
$$\bm x(\varepsilon)=\bm x_p(\varepsilon)+ \lambda_1(\varepsilon)\bm x_1(\varepsilon)+\ldots+\lambda_k(\varepsilon)\bm x_k(\varepsilon) $$

\noindent for some analytic functions $\lambda_1(\varepsilon),\ldots, \lambda_n(\varepsilon)$.

\end{theo}

\begin{theo} \label{theo.nohomoN}   Let $\bm\varepsilon_0\in \mathbb R^N$, for $N\geq 1$, and let $U$ be a neighborhood of this $\bm\varepsilon_0$.  Let $\bm D(\bm\varepsilon)=[\gamma_{ij}(\bm\varepsilon)]_{1\leq i,j\leq n}$, $\bm b(\bm\varepsilon)=(b_1(\bm\varepsilon),\ldots, b_n(\bm\varepsilon))^T$ be such that, for $1\leq i,j\leq n$, $\gamma_{ij}(\bm\varepsilon)$, $b_i(\bm\varepsilon)$ are in  $C^l(U)$, for $1\leq i,j\leq n$. Suppose that $r=\rank(\bm D(\bm\varepsilon_0))=\rank([\bm D(\bm\varepsilon_0)\mid \bm b(\bm\varepsilon_0)])$ and that for every $\bm\varepsilon\in U$,  $\rank([\bm D(\bm\varepsilon)\mid \bm b(\bm\varepsilon)])\leq r$. Let $k=n-r$.  Then there exist a neighborhood $U^*$ of $\bm\varepsilon_0$, a $C^l(U^*)$ vector field $\bm x_p(\bm\varepsilon)$ and a $C^l(U^*)$ frame field $(\bm x_1(\bm\varepsilon),\ldots, \bm x_k(\bm\varepsilon))$ such that any $C^l(U^*)$ vector field $\bm x(\bm\varepsilon)$ satisfying Equation \eqref{eq.1nh} can be writen as:
$$\bm x(\bm\varepsilon)=\bm x_p(\bm\varepsilon)+\lambda_1(\bm\varepsilon)\bm x_1(\bm\varepsilon)+\ldots+\lambda_k(\bm\varepsilon)\bm x_k(\bm\varepsilon) $$

\noindent for some functions $\lambda_1(\bm\varepsilon),\ldots, \lambda_n(\bm\varepsilon)$ in $C^l(U^*)$.

\end{theo}

The existence of $\bm x_p$ (the particular solution) is ensured by Theorem \ref{theo.nohomo} and the existence of the frame fields (homogeneous solutions) by Theorems \ref{theo.mainanalytic} and \ref{theo.main}.

\section{Sensitivity Analysis} \label{section.final}

In this section we study Problem \ref{problem.derivative}. The first two subsections deal with the case in which $\rank(\bm D(\bm \varepsilon_0))=n-1$ and  $\rank(\bm D(\bm \varepsilon))\leq n-1$ for $\bm \varepsilon$ in some neighborhood of $\bm \varepsilon_0$. In the first one, we present a direct method to solve Problem \ref{problem.derivative}. In the second one, an adjoint method is provided to perform this same task. This second type of methods are, computationally, more efficient for large values of $N$. Finally, in Subsection \ref{subsection.nosolution} we discuss the difficulties for computation of sensitivities in the case $\rank(\bm D(\bm \varepsilon_0))<n-1$.

\subsection{Direct Method} \label{subsection.si}
 Let us begin with the following straightforward observation. Recall that we are considering solutions verifying $\bm x(\bm \varepsilon_0)=\bm u$, for a given unitary vector $u$, as in the formulation of Problem \ref{problem.derivative}.

\begin{lemm} \label{lemm.equations} Let $\bm\varepsilon_0\in\mathbb R^N$ and $U$ be a open neighborhood of $\bm\varepsilon_0$. Let $\bm x(\bm\varepsilon)$ be a $C^1(\bm\varepsilon)$ vector field. Suppose that $\bm x(\bm\varepsilon)$ is a solution  of the system $\bm D(\bm\varepsilon)\bm x(\bm\varepsilon)=\bm b(\bm\varepsilon)$ for $\bm\varepsilon\in U$, where the entries in $\bm D(\bm\varepsilon)$ and $\bm b(\bm\varepsilon)$  are in  $C^1(U)$. Then:
\begin{equation} \label{eq.systemderivativemain} \bm D(\bm\varepsilon_0)\frac{\partial \bm x}{\partial \varepsilon_i}(\bm\varepsilon_0)=\frac{\partial \bm b}{\partial \varepsilon_i}(\bm\varepsilon_0)-\frac{\partial\bm D}{\partial \varepsilon_i}(\bm\varepsilon_0)\bm u.\end{equation}

\end{lemm}

\begin{proo} Equation \eqref{eq.systemderivativemain} is obtained taking derivatives from $\bm D(\bm\varepsilon)\bm x(\bm\varepsilon)=\bm b(\bm\varepsilon)$.

\hspace{12cm}\end{proo}

\begin{lemm} \label{lemm.equations} Let $\bm\varepsilon_0\in\mathbb R^N$ and let $U$ be a open neighborhood of $\bm\varepsilon_0$. Let $\bm x(\bm\varepsilon)$ be a $C^1(\bm\varepsilon)$ vector field. If  $\|\bm x(\bm \varepsilon)\|=1$ for $\bm\varepsilon\in U$, then, for $1\leq i\leq N$:
\begin{equation} \label{eq.unitaryderivative}\bm u^T \frac{\partial \bm x}{\partial \varepsilon_i}(\bm\varepsilon_0)=0.\end{equation}

\end{lemm}

\begin{proo}  Equation \eqref{eq.unitaryderivative} is obtained taking derivatives in $\bm x(\bm\varepsilon)^T  \bm x(\bm \varepsilon)=1$.

\hspace{12cm}\end{proo}

Combining Equations \eqref{eq.systemderivativemain} and \eqref{eq.unitaryderivative}   we can obtain a way to compute the derivatives inspired in Nelson's method (that was originally developed for eigenvector sensitivity, see 
\cite{B1, B2, N}).  Note that the general solution of Equation \eqref{eq.systemderivativemain} must be of the form
\begin{equation} \label{eq.systemderivativemainsolution}\frac{\partial\bm x}{\partial \varepsilon_i}(\bm\varepsilon_0)=\bm v+c \bm u  \end{equation}

\noindent for $\bm v$ being a particular solution of the system \eqref{eq.systemderivativemain} and $c$ a constant to be determined. From the fact that $\bm u^T\frac{\partial \bm x}{\partial \varepsilon_i}(\bm\varepsilon_0)=0$ we obtain that:
\begin{equation} \label{eq.systemc} c=-\bm v^T \bm u\end{equation}

\noindent So finally we get:

\begin{center} \begin{tabular}{| l l | } \hline & \textbf{Direct Method for the homogeneous system}\\
& $\bm D(\bm\varepsilon)\bm x(\bm\varepsilon)=\bm 0$ such that $\bm x(\bm\varepsilon_0)=\bm u$ and $\rank(\bm D(\varepsilon))=n-1$ is a neighborhood of $\varepsilon_0$.\\
\hline  1.- &We look for a particular solution $\bm v$ of the singular linear system $\bm D(\bm\varepsilon_0)\bm v=-\frac{\partial \bm D}{\partial \varepsilon_i}(\bm \varepsilon_0)\bm u$.\\
2.- &Set $c=-\bm v^T \bm u$.\\
3.- &Finally, $ \frac{\partial \bm x}{\partial \varepsilon_i}(\bm\varepsilon_0)=\bm v+c \bm u$.\\
\hline
\end{tabular}
\end{center}

\subsection{Adjoint Method}

Suppose that, for some $ F:\mathbb R^n\to \mathbb R$, for some $1\leq j\leq n$, we want to compute $\frac{\partial}{\partial \varepsilon_j}(F(\bm x(\bm\varepsilon_0)))$ where $\bm x(\bm \varepsilon)$ satisfies Equations \eqref{eq.1} and \eqref{eq.2}. Using this function $F$ we obtain some notational advantages, we study (effortlessly) a more general problem and we can easily particularize this problem (taking $F(\bm x(\bm \varepsilon))=x_i$) to recover $\frac{\partial  x_i}{\partial \varepsilon_j}(\bm \varepsilon_0)$ for $i=1,\ldots, n$.

For some  $\bm p\in\mathbb R^n$, $\lambda\in\mathbb R$ that do not depend on the variables in $\bm \varepsilon$,  consider the trivial equation:
$$ F(\bm x(\bm \varepsilon))=F(\bm x(\bm \varepsilon))+\underbrace{\bm p^T\bm D(\bm \varepsilon)\bm x(\bm \varepsilon)+\frac{1}{2}\bm \lambda(\bm x^T(\bm \varepsilon)\bm x(\bm \varepsilon)-1)}_{=0}$$

\noindent Taking derivatives we get:
$$\frac{\partial}{\partial \varepsilon_j}(  F(\bm x(\bm\varepsilon_0)))=\nabla  F(\bm x(\bm \varepsilon_0))\frac{\partial \bm x}{\partial \varepsilon_j}(\bm \varepsilon_0)+\bm p^T\left[\frac{\partial \bm D}{\partial \varepsilon_j}(\bm \varepsilon_0)\bm x(\bm \varepsilon_0)+\bm D(\bm \varepsilon_0)\frac{\partial \bm x}{\partial \varepsilon_j}(\bm \varepsilon_0)\right]+\bm \lambda\left[\bm x^T(\bm \varepsilon)\frac{\partial \bm x}{\partial  \varepsilon_j}(\bm \varepsilon)\right]$$

\noindent and re-organizing:
$$\frac{\partial}{\partial \varepsilon_j}(  F(\bm x(\bm\varepsilon_0)))=\bm p^T\frac{\partial \bm D}{\partial \varepsilon_j}(\bm \varepsilon_0)\bm x(\bm \varepsilon_0)+\underbrace{\left[\nabla  F(\bm x(\bm \varepsilon_0))+\bm p^T\bm D(\bm \varepsilon_0)+\bm \lambda\bm x^T(\bm \varepsilon)\right]}_{(\star\star)}\frac{\partial \bm x}{\partial \varepsilon_j}(\bm \varepsilon_0)$$

\noindent To simplify the expression above, we can choose $\bm p^T$ and $\lambda$ in such a way  that $(\star\star)=0$. To achieve this, we choose $\lambda\in\mathbb R$ to be the (unique) real number such that  the following linear system has a solution (at this point is where we need that $\rank(\bm D(\varepsilon))=n-1$) and then we solve it to find $\bm p$.
\begin{equation} \label{eq.morethanone} \bm D(\bm \varepsilon_0)^T\bm p=-\lambda\bm x(\varepsilon_0)-\nabla \bm F(\bm x(\bm \varepsilon_0))^T \end{equation}

Putting these ideas together we obtain:

\begin{center} \begin{tabular}{| l l | } \hline & \textbf{Adjoint Method for the homogeneous system}\\
& to compute partial derivatives of $F(\bm x(\bm\varepsilon))$ at $\bm \varepsilon_0$, where $\bm D(\bm\varepsilon)\bm x(\bm\varepsilon)=\bm 0$ and $\bm x(\bm\varepsilon_0)=\bm u$. \\
\hline  1.- & Find $\lambda$, $\bm p$ such that: \\
&$\bm D(\bm \varepsilon_0)^T\bm p=-\lambda\bm u-\nabla \bm F(\bm u)^T $.\\
2.- &$\frac{\partial}{\partial \varepsilon_j}( \bm F(\bm x(\bm \varepsilon_0)))=\bm p^T\frac{\partial \bm D}{\partial \varepsilon_j}(\bm \varepsilon_0)\bm u$.\\
\hline
\end{tabular}
\end{center}

Note that, $\bm p$ is not unique, that is, we can find two different solutions  $\bm p_1,\bm p_2$ satisfying \eqref{eq.morethanone} and so 
\begin{equation} \label{eq.pes} (\bm p_1^T-\bm p_2^T)\bm D=0\end{equation}

\noindent  But the value of $\frac{\partial}{\partial \varepsilon_j}( \bm F(\bm x(\bm \varepsilon_0)))$ does not vary. To check this, let us see that this invariance is equivalent to:
$$(\bm p_1^T-\bm p_2^T)\frac{\partial \bm D}{\partial \varepsilon_j}\bm x=0$$

\noindent and this equation is true as a consequence of Equations   \eqref{eq.systemderivativemain} and \eqref{eq.pes}.

If we want to compute several partial derivatives  $\frac{\partial \bm x}{\partial \varepsilon_{j_1}}(\bm \varepsilon_0),\ldots,\frac{\partial \bm x}{\partial \varepsilon_{j_r}}(\bm \varepsilon) $ the first part of the method is common to all of them and the second one just requires the computation of the corresponding derivative of $\bm D$ and a multiplication.

\subsection{Cases in which the solution of Problem \ref{problem.derivative} is not determined} \label{subsection.nosolution}

If for all $\varepsilon\in U$, $\rank(\bm D(\varepsilon))=n-k$, for $k>1$, then  Problem \ref{problem.derivative}  cannot be solved unless more information is provided. Let us see the following clarifying example:

\begin{exam} For the case $N=1, n=4$, $U=\mathbb R$ and $\varepsilon_0=0$. Consider  the matrix 
$$D(\varepsilon)=\begin{bmatrix}  \varepsilon & 0 & 0 & 0 \\ 0 & 0 & 0 & 0 \\ 0 & 0 & 0 & 0 \\ 0 & 0 & 0 & 0  \end{bmatrix}.$$
Consider some vector field $\bm x(\varepsilon)$ such that, for all $\varepsilon\in \mathbb R$, 
$$ \bm D(\varepsilon)\bm x(\varepsilon)=\bm 0,\qquad \bm x(0)=(0,1,0, 0)^T$$

\noindent This information is not sufficient to determine $\bm x'(0)$. To check this, just see that for every vector $(0,0,a,b)^T$ in the linear space spanned by $\{(0,0,1,0)^T,(0,0,0,1)^T\}$, the solution $\bm x(\varepsilon)=(0,\cos\varepsilon, \sin a\varepsilon, \sin b\varepsilon)$ satisfies $\bm x'(\varepsilon_0)=(0,0,a,b)$.

\end{exam}

\medskip

In the case of eigenproblems {this feature is well known}. It corresponds to the case in which $\lambda(\bm \varepsilon)$ is an eigenvalue of constant multiplicity $h>1$.

\medskip

Suppose that we have guaranteed the existence of a frame field $(\bm x_1(\bm \varepsilon),\ldots,\bm x_k(\bm\varepsilon))$ consisting of $k$ vector fields which are solutions of \eqref{eq.1} and that  some orthonormal vectors $\bm u_1,\ldots, \bm u_k$ are provided, so that it is prescribed that $\bm x_j(\bm \varepsilon_0)=u_j$ for $j=1\ldots, k$. In this case, the derivative of the solutions must satisfy some extra conditions. For $1\leq j_1,j_2\leq k$, $j_1\neq j_2$:
$$\bm x_{j_1}(\bm \varepsilon)\bm x_{j_2}(\bm \varepsilon)=0\Rightarrow \bm u_{j_1}^T\frac{\partial\bm x }{\partial\varepsilon_i }(\bm \varepsilon_0)+\bm u_{j_2}^T\frac{\partial \bm x }{\partial \varepsilon_i}(\bm \varepsilon_0)=0  $$

\noindent Even with this extra conditions, the solution to Problem \eqref{problem.derivative} is not determined. It can be proved in a straightforward manner, studying the corresponding linear system, that has $\binom{k}{2}$ degrees of freedom, or just noting the following:

\begin{rem} In the notation of Remark \ref{rem.movimiento}, see that for the frame field $(\bm x_1(\bm \varepsilon),\ldots,\bm x_k(\bm \varepsilon))$ obtained following the method in the proof of Theorems \ref{theo.mainanalytic} and  \ref{theo.main} any other frame field $(\bm y_1(\bm\varepsilon),\ldots, \bm y_k(\bm \varepsilon))$ satisfies the corresponding conditions  if and only if 
\begin{equation} \label{eq.lie1} \text{for } 1\leq i\leq j, \;\bm y_i(\varepsilon)=\bm K(\bm\varepsilon)\bm x_i(\bm\varepsilon)\end{equation}

Now, if we impose that, for $1\leq i\leq k$, $\bm x_i(\bm\varepsilon_0)=\bm y_i(\bm\varepsilon_0)=\bm u_i$ for some prescribed vector $\bm u_i$, then $\bm K(\bm\varepsilon_0)=\bm I_{n}$ and so $\bm K(\bm\varepsilon)$ is in $\mathcal{SO}(n)$. Taking derivatives in Equation \eqref{eq.lie1} we have that
$$\frac{\partial\bm y_i}{\partial \varepsilon_i}(\varepsilon_0)=\frac{\partial \bm K}{\partial \varepsilon_i}(\bm\varepsilon_0)\bm u_i+\frac{\partial\bm x_i}{\partial \varepsilon_i}(\bm\varepsilon)$$

\noindent If we want $\bm K(\bm\varepsilon)$ to be of the type explained in Equation \eqref{eq.K} and to satisfy $\bm K(\bm\varepsilon_0)=\bm I_{n}$ we still have $\binom{k}{2}$ free parameters since $\frac{\partial \bm T}{\partial \varepsilon_i}(\varepsilon_0)$ can be any skew-symmetric matrix of size $k\times k$ (see \cite{BAKER}).

\end{rem}

\section*{Ackwnoledgements}

This work has been supported by the Spanish Agencia Estatal de Investigaci\'on through project PID2020-116207GB-I00  and Junta de Comunidades de Castilla-La Mancha through project  SBPLY/19/180501/000110.

\section*{References}


\end{document}